\documentclass[12pt]{article}
\usepackage{amssymb, amscd, amsmath, amsthm, epsf, epsfig, latexsym}

\def\zz{{\bf Z}}
\def\qq{{\bf Q}}

\newtheorem{thm}{Theorem}[section]

\newtheorem{defn}[thm]{Definition}

\title{Pseudo--Slice Knots}

\author{C. Livingston\thanks{Dept. of Mathematics,
Indiana University, Bloomington, IN 47405, 
livingst@indiana.edu}  }

\begin{document}
\maketitle
\begin{abstract}

For $n >1$, if the Seifert form of a knotted $2n-1$--sphere $K$ in $S^{2n+1}$ has a
metabolizer, then the knot is slice.  Casson and Gordon proved that this is false in
dimension three.  However, in the three dimensional case it is true that if the
metabolizer has a basis represented by a strongly slice link then $K$ is slice. The
question has been asked as to whether it is sufficient that each basis element is
represented by a slice knot to assure that $K$ is slice.  For genus one knots this is of
course true; here we present a genus two counterexample. 
   
\end{abstract}


\section{Introduction.} Let $K$ be a knot in $S^3$. If $F$ is a Seifert surface for $K$
there is a Seifert form,
$V_K$, defined on
$H_1(F,\zz)$. The knot $K$ is called {\sl algebraically slice} if $V_K$ 
vanishes on some half--dimensional summand of $H_1(F,\zz)$; such a summand is called a
{\sl metabolizer} for
$V_K$.   If $K$ is {\sl slice},
that is, if it bounds a smooth embedded disk in $B^4$, then it is algebraically slice. 
Casson and Gordon \cite{cg1,cg2} proved that the converse does not hold by constructing
explicit examples of algebraically slice knots that are not slice.  (This was in contrast
to the result proved by Levine \cite{le1, le2} that in higher dimensions the analagous
condition of algebraic sliceness does imply that a knot is slice.)

If a basis for the metabolizer in $H_1(F,\zz)$ is represented by a stongly slice link, one
for which the components bound disjoint disks in $B^4$, then it is easily shown that $K$ is
slice.  The question has been asked whether it is sufficient to show that a basis of the
metabolizer is represented by slice knots to assure that $K$ will be slice.  Here we
provide a genus two counterexample to show this is not the case.  Litherland has
previously given such an example in
\cite{lit1}, but that paper never appeared, and the result depended on the development of
a lengthy algorithm for computing Casson--Gordon invariants; that algorithm itself has not
appeared in print.

The example here points to a much deeper question regarding classical concordance.  At an
empirical level, knotting in the curves representing elements in the metabolizer of the
Seifert form present secondary obstructions to slicing a classical knot.  This was first
made formal in Gilmer's work \cite{gi} where certain signatures of these knots were
related to Casson-Gordon invariants.  In the case we are considering these signatures all
vanish, and hence a more subtle approach is needed.  It is expected that in addition to
the signatures of the individual components of the metabolizing basis there should be
abelian invariants of the entire link that provide second order slicing obstructions.  The
example produced here points to the existance of such obstructions, but the precise
formulation is not evident and remains an open question.

The results of this paper depend only on the original definition of Casson--Gordon
invariants of \cite{cg1}, the connected sum formula of Gilmer \cite{gi}, and a simple
method for computing Casson--Gordon invariants of satellite knots,
first described in \cite{lit2} and reformulated in \cite{gl1} as needed here.

See \cite{rol} and \cite{bz} for basic results in knot theory.  We will work in the smooth
category, but all results extend to the topological, locally flat, category by \cite{fq}.

The example presented here was developed in response to a question of Effie Kalfagianni
addressed to the author.  We also wish to thank Pat Gilmer for conversations regarding
this work.


\section{Casson--Gordon Invariants.} 

Let $K$ be a knot in $S^3$ with 2--fold branched cover $M_K$.  For a character
$\chi:H_1(M_K,\zz) \to \zz_p$, Casson and Gordon \cite{cg1} define an
invariant, denoted $\sigma_1\tau(K,\chi) \in \qq$.  To simplify notation we have:

\begin{defn}  $\sigma(K,\chi) =  \sigma_1\tau(K,\chi)$. \end{defn}

There is a linking form lk:$H_1(M_K,\zz) \times H_1(M_K,\zz) \rightarrow \qq/\zz$.
The main result of \cite{cg1} is the following.

\begin{thm}  If $K$ is slice there is a subgroup $H \subset H_1(M_K,\zz)$ with
$|H|^2 = |H_1(M_K,\zz)|$ and such that for any character $\chi$ with values in $\zz_p$
for prime $p$ that vanishes on $H$, $\sigma(K,\chi) = 0$.  Furthermore $H$ can be assumed
to be a metabolizer (self--annihilating) for the linking form on $M$.

\end{thm}

One simple result concerning the Casson--Gordon invariant is that $\sigma(K,\chi) = 
\sigma(K,-\chi)$.  A much deeper result is the additivity result proved by Gilmer
\cite{gi}, as we now describe.  If
$K = J_1 \# J_2$, then $H_1(M_K,\zz) =  H_1(M_{J_1},\zz) \oplus H_1(M_{J_2},\zz)$ and any
character $\chi$ on $ H_1(M_{K},\zz)$ can be written as $\chi = \chi_1 \oplus \chi_2$ with
$\chi_i $ a character on $H_1(M_{J_i},\zz)$.

\begin{thm} $\sigma(K,\chi) = \sigma(J_1,\chi_1) +
\sigma(J_2,\chi_2)$. 
\end{thm}

For satellite knots there is an algorithm that simplifies the computation of
its Casson--Gordon invariants. Details appear in \cite{gl1, lit2}; here is a summary. 
Suppose that
$L$ is an unknotted circle in
$S^3$ in the complement of $K$ that is null homologous in $S^3 - K$.  If a neighborood of
$L$ is removed from $S^3$ and replaced with the complement of a knotted circle $J$ in
$S^3$ (with the boundaries identified so that the longitude of $L$ is identified with the
meridian of $J$ and the meridian of $L$ is identified with the longitude of $J$) then the
resulting manifold is diffeomorphic to
$S^3$, but the curve
$K$ now represents a perhaps different knot, say
$K^*$, in $S^3$.  (Traditionally $K^*$ has been called a satellite of $J$ with
embellishment $K$; in effect the portion of $K$ that passes through $L$ is tied into the
knot $J$.)

The curve $L$ lifts to a pair of curves, $\tilde{L'}$ and  $\tilde{L''}$ in $M_K$.  Thus,
$M_{K^*}$ is constructed from $M_{K}$ by removing neighborhoods of  $\tilde{L'}$ and 
$\tilde{L''}$ and replacing both with copies of the complement of $J$.  This construction
leaves the homology unchanged and there is a natural correspondence between the homology,
and cohomology, groups of $M_{K}$ and $M_{K^*}$. In particular we can identify a characters
$\chi$ on $H_1(M_{K^*},\zz)$ with characters on $H_1(M_{K},\zz)$.  In this
situation we have the following:

\begin{thm}\label{surgeryformula} $\sigma(K^*,\chi) = \sigma(K,\chi) + 2\sigma_{
\chi(\tilde{L'})/p}(J)$.
\end{thm}

Here $\sigma_{k/p}(J)$ denotes the classical Tristram-Levine signature  \cite{tr} of
the knot $J$, given as the signature of the hermetianized Seifert form $(1 - \omega)V_J +
(1- \bar{\omega})V_J^t$, $\omega = e^{ k2\pi i /p}$.

A similar result holds if $L$ is replaced with a nullhomologous unlink in the complement
of $K$, with the single signature in Theorem \ref{surgeryformula} replaced with a sum of
signatures.

\section{A Pseudo-Slicable Genus Two Knot.}
In this section we construct an example of a genus two knot, $K^*$, that
is algebraically slice and for which a generating set for a metabolizer of a Seifert
form can be represented by slice knots, and yet the knot itself is not slice.


\vskip.2in
\epsfxsize=4in
\centerline{\epsfbox{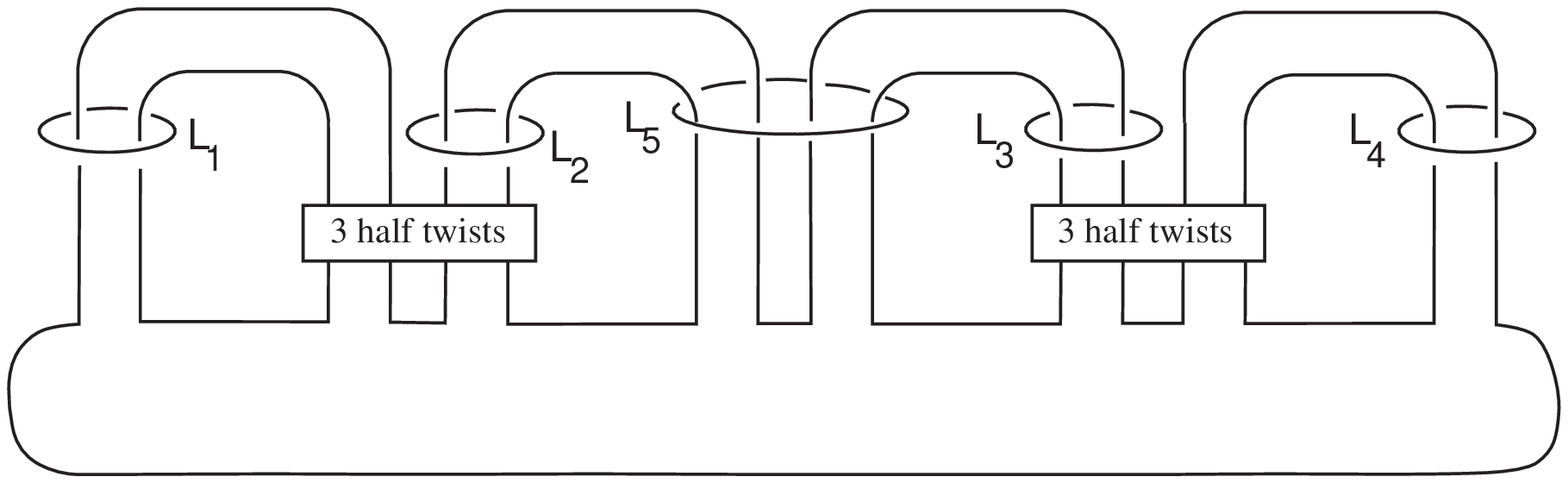}}
\vskip.1in
\centerline{Figure 1. }
\vskip.2in

The example is illustrated in Figure 1.  In  the figure,  a knot $K$
is drawn along with five curves in the complement of its Seifert
surface, 
$L_i, i = 1,
\ldots 5$.  The bands are twisted in such a
way that the Seifert form for this knot is given by
$$\left( \begin{matrix}  0 & 1 \\ 2 & 0 
\end{matrix} \right) \oplus \left( \begin{matrix}  0 & 1 \\ 2 & 0 
\end{matrix} \right).$$
The knot $K^*$ is constructed by removing neighborhoods of the $L_i$ and
replacing them with knot complements.  In our case all of these will be the
complement of the same knot, $J$, with the exception of $L_5$ which is
replaced with the complement of $-J$. 

Let $\{x_1, y_1, x_2, y_2\}$ denote the standard symplectic basis for the
first homology of the Seifert surface for $K^*$ (as drawn above) and at
the same time denote obvious simple curves on the Seifert surface
representing these classes.  (So, each curve passes over exactly one
band, with $x_1$ going around the leftmost band.)

Notice that, as knots, $x_1$ and $y_2$ are each the knot $J$, while $y_1$
and $x_2$ are each represented by the slice knot $J \# -J$.  In
particular, $y_1$ and $x_2$ form a basis for a metabolizer for the
Seifert form with both elements represented by slice knots.  That this
 knot provides the desired examples follows from the following theorem.

\begin{thm} There exists a number $C >0 $ such that if $\sigma_{1/3}(J)
> C/2
$, then
$K^*$ is not slice.

\end{thm}

To prove this result, we begin with a diagram of the 2-fold branched cover
of $S^3$ branched over $K^*$ as drawn in Figure 2 using the algorithm of
\cite{ak}.  In this figure the lifts of the $L_i$ are drawn, and are to
be replaced with the appropriate knot complement (of either $J$ or $-J$)
to complete the construction.


\vskip.2in
\epsfxsize=4.5in
\centerline{\epsfbox{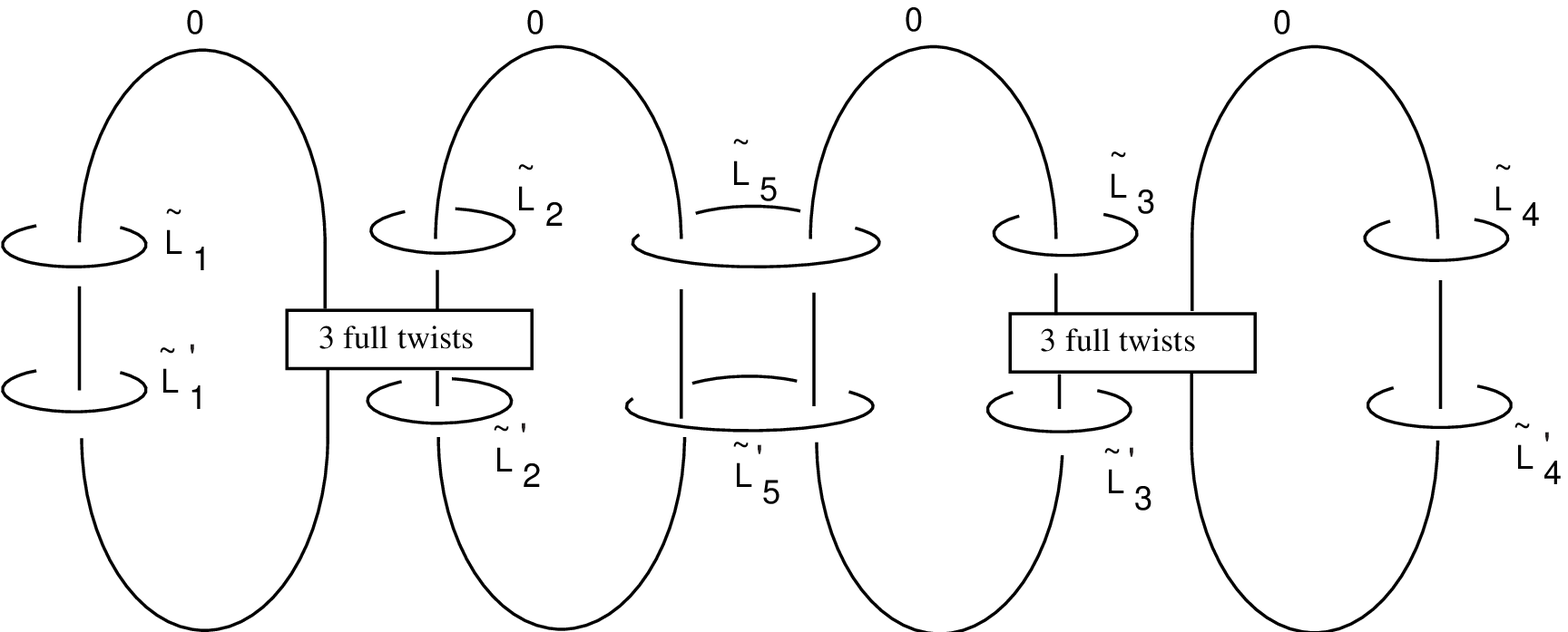}}
\vskip.1in
\centerline{Figure 2. }
\vskip.2in

The homology of the cover is $(\zz_3)^4$ generated by meridians to the
four surgery circles.  Denote by $\chi_{(a,b,c,d)} $ the $\zz_3$ valued
character taking values $a, b, c$, and $d$ on the four meridians,
listed as drawn from left to right, respectively.  By choosing the
orientations of the $\tilde{L}$ properly, we can assume that:

\begin{itemize}
\item $\chi_{(1,0,0,0)}$ takes value 1 on $\tilde{L_1}$ and value 0 on all
other $\tilde{L_i}$.
\item $\chi_{(0,1,0,0)}$ takes value 1 on $\tilde{L_2}$ and $\tilde{L_5}$
and value 0 on all other $\tilde{L_i}$.
\item $\chi_{(0,0,1,0)}$ takes value 1 on $\tilde{L_3}$ and $\tilde{L_5}$
and value 0 on all other $\tilde{L_i}$.
\item $\chi_{(0,0,0,1)}$ takes value 1 on $\tilde{L_4}$ and value 0 on all
other $\tilde{L_i}$.
\end{itemize}
The same values are taken on the translates of the $\tilde{L_i}$,
denoted by $\tilde{L'_i}$ in the figure. 

From this calculation we have the following formula:
$$\sigma(K^*,\chi_{(a,b,c,d)}) = \sigma(K,\chi_{(a,b,c,d)}) +
2(\overline{a} +
\overline{b} + \overline{c} + \overline{d} -
\overline{b+c})\sigma_{1/3}(J)$$
where $\overline{x} = 1$ if $x = \pm 1 \in \zz_3$ and $\overline{x} = 0$
if $x = 0$.

Since the set of values of $\{ \sigma(K,\chi_{(a,b,c,d)})\}_{(a,b,c,d) \in (\zz_3)^4}$
represents a fixed finite set of rational numbers, it is bounded above in absolute value
by some number
$C$. It is clear that if we can prove that for every possible metabolizer
the value of
$\overline{a} +
\overline{b} + \overline{c} + \overline{d} -
\overline{b+c}$ is positive for some character vanishing on the
metabolizer, then if $\sigma_{1/3}(J) > C/2$ the corresponding
Casson--Gordon invariant is nonzero.

It is easily seen that for all $b$ and $c$, $\overline{b} + \overline{c}  -
\overline{b+c} \ge 0$.  Hence, $\overline{a} +
\overline{b} + \overline{c} + \overline{d} -
\overline{b+c}$ will be positive unless both $a$ and $d$ are 0.  But if
this were the case for all characters vanishing on the metabolizer,
then this set of characters would necessarily be spanned by the
characters $\chi_{(0,1,0,0)}$ and $\chi_{(0,0,1,0)}$.  But then the set
of characters would also contain $\chi_{(0,1,1,0)}$, and for this
character $\overline{a} +
\overline{b} + \overline{c} + \overline{d} -
\overline{b+c} = 1$.  

We have now seen that for every metabolizer there is some character
vanishing on that metabolizer for which $\overline{a} +
\overline{b} + \overline{c} + \overline{d} -
\overline{b+c} >0$.  Since $\sigma(K,\chi_{(a,b,c,d)}) \ge -C$ and
$\sigma_{1/3}(J) > C/2$, it certainly follows that for this character
$\sigma(K^*,\chi_{(a,b,c,d)}) > 0$ and the proof is complete.

\vskip.5in


\newcommand{\etalchar}[1]{$^{#1}$}

\end{document}